\documentclass[11pt]{amsart}

\usepackage[T1]{fontenc}
\usepackage[utf8]{inputenc}
\usepackage{amsmath,amssymb,mathtools}
\usepackage{lmodern}
\usepackage[margin=1in]{geometry}
\usepackage[final]{microtype}
\usepackage{needspace}
\usepackage{hyperref}
\usepackage{bookmark}

\hypersetup{
  hidelinks,
  bookmarksnumbered=true,
  pdfdisplaydoctitle=true,
  pdftitle={Joint Extension Field Point Counts for Superelliptic Curves},
  pdfauthor={Yuxiang Yao and Xu Zhuang}
}
\renewcommand{\MR}[1]{}

\numberwithin{equation}{section}
\allowdisplaybreaks[3]

\newtheorem{theorem}{Theorem}[section]
\newtheorem{proposition}[theorem]{Proposition}
\newtheorem{corollary}[theorem]{Corollary}
\newtheorem{lemma}[theorem]{Lemma}
\theoremstyle{definition}
\newtheorem{definition}[theorem]{Definition}
\newtheorem{example}[theorem]{Example}
\theoremstyle{remark}
\newtheorem{remark}[theorem]{Remark}

\newcommand{\A}{\mathbf A}
\newcommand{\C}{\mathbf C}
\newcommand{\F}{\mathbf F}

\newcommand{\N}{\mathbf N}
\newcommand{\PP}{\mathbf P}

\newcommand{\Z}{\mathbf Z}

\newcommand{\E}{\mathbb E}
\newcommand{\Prb}{\mathbb P}
\newcommand{\cO}{\mathcal O}
\newcommand{\frakm}{\mathfrak m}
\DeclareMathOperator{\ev}{ev}
\newcommand{\Exp}{\operatorname{Exp}_{\sigma}}
\newcommand{\Log}{\operatorname{Log}_{\sigma}}
\DeclareMathOperator{\Cov}{Cov}
\DeclareMathOperator{\Var}{Var}
\DeclareMathOperator{\Fil}{Fil}

\DeclareMathOperator{\lcmop}{lcm}

\DeclareMathOperator{\gh}{gh}
\newcommand{\aff}{\operatorname{aff}}

\title[Joint extension field point counts]{%
  Joint Extension Field Point Counts for Superelliptic Curves}
\author{Yuxiang Yao}
\address{Department of Mathematics\newline
  \indent University of California, Irvine\newline
  \indent Irvine, CA 92697-3875, USA}
\email{yuxiany@uci.edu}
\author{Xu Zhuang}
\address{Department of Mathematics\newline
  \indent University of California, Irvine\newline
  \indent Irvine, CA 92697-3875, USA}
\email{xzhuang8@uci.edu}
\date{}

\begin{document}

\begin{abstract}
Fix a finite field \(\F_q\), an integer \(m\ge2\) prime to \(q\), and
\(2\le\kappa\le m\). For a uniformly chosen \(\kappa\)-power-free
polynomial \(F\in\F_q[x]\) of degree \(D\), we determine the limiting
joint distribution of the affine point counts of \(y^m=F(x)\) over
\(\F_q,\ldots,\F_{q^n}\) as \(D\to\infty\), with \(q\) and \(n\) fixed.
The probability generating function is a finite Euler product with explicit
local factors. Equivalently, the limiting
vector is a sum of independent closed-point contributions, each of which
may affect several extension degrees. We obtain the formula by specializing
the equidistribution and motivic Euler product formalism of Bertucci--Howe.
We then compute means, covariances, and moment expansions as \(q\) varies
in a fixed reduced residue class modulo \(m\), after the degree limit has
been taken. The one-field specialization recovers the power-free affine
point-count distribution of Cheong--Wood--Zaman.
\end{abstract}

\maketitle

\section{Introduction}\label{sec:introduction}

Let \(q\) be a prime power, let \(m\ge2\) be prime to \(q\), and fix an
integer \(\kappa\) with \(2\le \kappa\le m\).  For \(D\ge1\), let \(U_D\) be the
set of \(\kappa\)-power-free polynomials \(F\in\F_q[x]\) of degree exactly
\(D\), equipped with the uniform probability measure.  Thus \(F\in U_D\) if
\(\deg F=D\) and no irreducible polynomial \(P\in\F_q[x]\) satisfies
\(P^\kappa\mid F\).  For \(F\in U_D\), consider the affine superelliptic model
\[
  C_F^{\aff}:\quad y^m=F(x).
\]
For \(r\ge1\), write
\[
  N_r(F):=\#C_F^{\aff}(\F_{q^r})
  =
  \#\{(a,y)\in\F_{q^r}^2:y^m=F(a)\}.
\]
We study, for fixed \(n\), the limiting joint distribution of
\((N_1(F),\ldots,N_n(F))\) as \(D\to\infty\).

The polynomials are not required to be monic, and the counts refer to the
affine models as written, including their singular points. We impose no
geometric irreducibility condition. Throughout the degree limit,
\(q,m,\kappa\), and \(n\) remain fixed.

The distribution of rational points in families over a fixed finite field
has been studied in several settings. Kurlberg and Rudnick~\cite{KR}
showed that the large-genus distribution in the squarefree hyperelliptic
family is governed by independent local random variables. Bucur, David,
Feigon, and Lal\'in~\cite{BDFL} established analogous results for cyclic
trigonal curves and discussed prime-degree cyclic covers. For a prime
\(\ell\) with \(q\equiv1\pmod\ell\), Bucur, David, Feigon, Kaplan,
Lal\'in, Ozman, and Wood~\cite{BDFKLOW} determined the point-count
distribution over the full moduli space of cyclic \(\ell\)-covers of a
given genus. Meisner~\cite{Meisner} treated more general cyclic covers by
letting the degrees of all squarefree factors in their defining polynomials
tend to infinity. 

The closest one-field comparison is Cheong--Wood--Zaman~\cite{CWZ}.
Their Theorem~3.1 determines the limiting distribution of points on
power-free affine superelliptic models and on their normalizations as the
polynomial degree grows. Their power-free exponent corresponds to
\(\kappa\) here. Thus the \(n=1\) specialization of our result recovers
their affine distribution for every \(2\le\kappa\le m\). In the related monic squarefree family,
Chen~\cite{Chen} used the homology of braid groups and the Burau
representation to obtain an exact first-moment identity.

Point counts over extensions of the ground field have also been considered.
Alzahrani~\cite{Alzahrani} studied the distribution over a fixed finite
extension in the hyperelliptic family. Here we keep several extension
degrees simultaneously. Their counts are linked because a closed point of
degree \(e\) contributes to every \(N_r\) with \(e\mid r\), and all these
contributions depend on the same local residue. The joint law therefore
requires the simultaneous power-residue conditions in the fields
\(\F_{q^r}\). Our local formula records these conditions through subgroup
indices in \(\F_{q^e}^{\times}\) and inclusion--exclusion.

Our proof builds on Howe's probabilistic formalism for random variables in
\(\lambda\)-rings~\cite[Sections~3 and 5]{Howe2024}.
In this formalism, the \(\Lambda\)-distribution of \(X\) records the
expectations \(\E[f\circ X]\), where \(f\) ranges over the ring of symmetric
functions. These expectations are encoded by the \(\sigma\)-moment
generating function, defined using the plethystic exponential.
The falling moment generating function encodes the same
\(\Lambda\)-distribution~\cite[Example~3.2.3]{Howe2024}.

Bertucci--Howe~\cite[Theorem~4.2.1]{BH} show that, under an
equidistribution hypothesis, the limiting \(\sigma\)-moment generating
function is a motivic Euler product. We use the joint falling moment
version of this result~\cite[Remark~4.2.2]{BH}.
In the present family, finitely many
nonzero local \(\kappa\)-jets become independent and uniform as
\(D\to\infty\). We apply the first ghost coordinate and a continuous
specialization of the symmetric-function variables to obtain the ordinary
probability generating function of \((N_1,\ldots,N_n)\). Only closed
points of degree at most \(n\) contribute, so the resulting Euler product
is finite.

The block description gives exact means and covariances, as well as moment
expansions in \(q^{-1}\). These expansions are taken only after the
\(D\to\infty\) limit, with \(q\) restricted to a reduced residue class
modulo \(m\). In the hyperelliptic example, the mean over an even-degree
extension has a leading excess of size \(q^{r/2}\). This is related to the
phenomenon investigated by Brock and Granville~\cite{BrockGranville} and
explained by Katz~\cite{KatzQuadraticExcess} through monodromy and
Frobenius--Schur indicators. Our calculation concerns affine models in the
degree-first limit, so we compare the excess terms while keeping the
families and the orders of limits explicit. See
Example~\ref{ex:squarefree-hyperelliptic} and the remark following it.

We now introduce the notation used in the main formula.  Let \(a_e\) be the
number of degree-\(e\) closed points of \(\A^1_{\F_q}\).  Equivalently,
\[
  a_e=\frac1e\sum_{\ell\mid e}\mu(\ell)q^{e/\ell}.
\]
For fixed \(n\), put
\[
  R_e=R_e(n):=\{1\le r\le n:e\mid r\}.
\]
For \(e\mid r\), define
\[
  \begin{aligned}
    g_r&:=\gcd(m,q^r-1),
    & c_r^{(e)}&:=\frac{q^r-1}{q^e-1},\\[4pt]
    \gamma_r^{(e)}&:=\gcd(g_r,c_r^{(e)}),
    & \delta_r^{(e)}&:=\frac{g_r}{\gamma_r^{(e)}}.
  \end{aligned}
\]
Equivalently, if
\[
  H_r^{(e)}:=(\F_{q^r}^{\times})^m\cap\F_{q^e}^{\times},
\]
then \(\delta_r^{(e)}\) is the index of \(H_r^{(e)}\) in
\(\F_{q^e}^{\times}\), and \(\gamma_r^{(e)}=g_r/\delta_r^{(e)}\).  The number
\(\gamma_r^{(e)}\) is the expected number of nonzero solutions to \(y^m=u\) over
\(\F_{q^r}\), when \(u\) is chosen uniformly from \(\F_{q^e}^{\times}\).

For a degree-\(e\) closed point, set
\[
  \alpha_e:=\frac{q^{e(\kappa-1)}-1}{q^{e\kappa}-1}.
\]
Here \(\alpha_e\) is the probability that a uniformly chosen nonzero
\(\kappa\)-jet over \(\F_{q^e}\) has positive valuation.  Finally, for
\(S\subseteq R_e\), let \(C_S^{(e)}\) be the number of units
\(u\in\F_{q^e}^{\times}\) whose membership pattern is exactly \(S\), meaning
\[
  u\in H_r^{(e)}
  \quad\Longleftrightarrow\quad
  r\in S.
\]

We take an empty product to be \(1\). The main theorem is the following.

\begin{theorem}\label{thm:main-pgf}
For fixed \(n\), with respect to the ordinary uniform probability measure on
\(U_D\),
\[
  \lim_{D\to\infty}
  \E_{F\in U_D}\prod_{r=1}^n w_r^{N_r(F)}
  =
  G_n(\mathbf w):=\prod_{e=1}^n L_e(\mathbf w)^{a_e},
\]
where \(\mathbf w=(w_1,\ldots,w_n)\) and
\[
  L_e(\mathbf w)
  =\alpha_e\prod_{r\in R_e}w_r^e
  +(1-\alpha_e)\sum_{S\subseteq R_e}
    \frac{C_S^{(e)}}{q^e-1}\prod_{r\in S}w_r^{eg_r}.
\]
Equivalently, if
\[
  \mathbf N_D(F):=(N_1(F),\ldots,N_n(F)),
\]
then \(\mathbf N_D\Rightarrow\mathbf N^*\), where
\(\mathbf N^*=(N_1^*,\ldots,N_n^*)\) is the
\(\Z_{\ge0}^n\)-valued random vector with probability generating function
\(G_n\).  In particular, for each
\(\mathbf k=(k_1,\ldots,k_n)\),
\[
  \lim_{D\to\infty}
  \Prb_{F\in U_D}\bigl(\mathbf N_D(F)=\mathbf k\bigr)
  =
  [w_1^{k_1}\cdots w_n^{k_n}]\,G_n(\mathbf w).
\]
\end{theorem}

Since \(0\le N_r(F)\le g_rq^r\), the vectors
\((N_1(F),\ldots,N_n(F))\) all lie in the same finite set.  Thus coefficientwise
convergence of the probability generating functions is equivalent to convergence
of probability masses, and convergence in distribution also implies convergence of
all fixed moments.

The limiting distribution admits a useful independent-block description.  For
\(e\le n\), let
\[
  X_e=(X_e^{(1)},\ldots,X_e^{(n)})
\]
be the \(\Z_{\ge0}^n\)-valued random vector defined as follows.  Take independent
random variables \(\eta_e\) and \(U_e\), where
\[
  \Prb(\eta_e=1)=\alpha_e,\qquad
  \Prb(\eta_e=0)=1-\alpha_e,
  \qquad
  U_e\sim\operatorname{Unif}(\F_{q^e}^{\times}).
\]
Then set
\[
  X_e^{(r)}
  =
  \begin{cases}
    e\eta_e+
    eg_r(1-\eta_e)\mathbf 1_{\{U_e\in H_r^{(e)}\}},
    & e\mid r,\\
    0, & e\nmid r.
  \end{cases}
\]
The superscript \((r)\) denotes the \(r\)-th coordinate.  This random vector has
probability generating function \(L_e(\mathbf w)\).

\begin{corollary}\label{cor:blocks}
For each \(1\le e\le n\), let \(X_{e,1},\ldots,X_{e,a_e}\) be independent
copies of \(X_e\), and take all these blocks to be mutually independent as
\(e\) varies.  Then
\[
  \mathbf N^*
  \stackrel{\mathrm{law}}{=}
  \sum_{e=1}^n\sum_{j=1}^{a_e}X_{e,j}.
\]
In particular, the limiting law is a sum of independent closed-point
contributions, and a degree-\(e\) block contributes only to the coordinates
\(N_r^*\) with \(e\mid r\).
\end{corollary}

Let us spell out the specialization \(n=1\).  Then only degree-\(1\) closed
points contribute: \(R_1=\{1\}\), \(a_1=q\), and
\(g_1=\gcd(m,q-1)\).  Thus \(X_1\) is a scalar random variable.  In the
squarefree case \(\kappa=2\), one has
\[
  \alpha_1=\frac{q^{-1}}{1+q^{-1}},
  \qquad
  1-\alpha_1=\frac1{1+q^{-1}}.
\]
Since \((\F_q^\times)^m\) has index \(g_1\) in \(\F_q^\times\), the distribution
of \(X_1\) is
\[
  X_1=
  \begin{cases}
    0, & \text{with probability }
    \frac1{1+q^{-1}}\left(1-\frac1{g_1}\right),\\[4pt]
    1, & \text{with probability }
    \frac{q^{-1}}{1+q^{-1}},\\[4pt]
    g_1, & \text{with probability }
    \frac1{(1+q^{-1})g_1}.
  \end{cases}
\]
When \(g_1=1\), the last two values coincide and the corresponding probabilities
are added.  Therefore
\[
  N_1^*
  \stackrel{\mathrm{law}}{=}
  \sum_{a\in\F_q}X_{1,a},
\]
where the \(X_{1,a}\) are independent copies of \(X_1\).  This recovers the
one-field limiting distribution of Cheong--Wood--Zaman in the squarefree
specialization, and Theorem~\ref{thm:main-pgf} extends it to the joint vector
\((N_1,\ldots,N_n)\). Corollary~\ref{cor:blocks} is also the main tool for the moment
calculations later in the paper.

The paper is organized as follows.  Section~\ref{sec:preliminaries} recalls the
Bertucci--Howe equidistribution and motivic Euler product formalism needed for
the proof.  Section~\ref{sec:local-factors} computes the local point-counting
factors.  Section~\ref{sec:main-proof} proves Theorem~\ref{thm:main-pgf} and
the independent-block corollary.  Section~\ref{sec:moments} applies the block
model to moments, covariances, and \(q\)-asymptotic expansions.

\section*{Acknowledgments}
We thank Sean Howe for suggesting that we apply \(\Lambda\)-distributions
and the equidistribution methods of Bertucci--Howe to point-count
statistics of superelliptic curves, and for his helpful advice on this
project. We also thank Nathan Kaplan for pointing out the
connection between the even degree excess in
Example~\ref{ex:squarefree-hyperelliptic} and the
Brock--Granville excess phenomenon~\cite{BrockGranville}.

\section{Preliminaries on Bertucci--Howe equidistribution}\label{sec:preliminaries}

We recall only the part of the preliminaries needed for the present
calculation. The main references for this section are Howe~\cite{Howe2024} and Bertucci--Howe~\cite{BH}.

\subsection{Admissible \texorpdfstring{\(\Z\)}{Z}-sets}

\begin{definition}\label{def:admissible-Z-set}
An \emph{admissible $\Z$-set} is a set $V$ with a $\Z$-action such that
$V=\bigcup_{k\ge1}V^{k\Z}$ and
$V^{k\Z}:=\{v\in V\mid k\cdot v=v\}$ is finite for every $k\ge1$.
We write $|V|$ for the set of $\Z$-orbits.  For $v\in V$, the degree
$\deg(v)=\deg(|v|)$ is the cardinality of its orbit.
\end{definition}

For $k\ge1$, let $V_k$ denote the same underlying set with the $\Z$-action
multiplied by $k$.  Thus $V_k(1)$, the set of fixed points of $V_k$, is
canonically $V^{k\Z}$.  If $\varphi:V\to B$ is a morphism of admissible
$\Z$-sets, then $\varphi_k:V_k\to B_k$ denotes the induced morphism.  For
$b\in B$, the fiber $V_b$ is equipped with the action obtained by multiplying
the action on $V$ by $\deg(b)$.

For example, let $V=Y(\overline{\F}_q)$ for a variety
$Y/\F_q$, with $1\in\Z$ acting by geometric Frobenius. Then $|V|$ is the set
of closed points of $Y$, degrees are residue-field degrees, and $V_k$ is the
admissible $\Z$-set associated with the base change $Y_{\F_{q^k}}$.

\subsection{\texorpdfstring{\(W(\C)\)}{W(C)}-valued functions and \texorpdfstring{\(\lambda\)}{lambda}-probability spaces}

Let $W(\C)$ be the big Witt ring.  We write
$\gh_k:W(\C)\to\C$ for its $k$th ghost coordinate and $p_j\circ(-)$ for its
$j$th Adams operation, where $p_j$ is the $j$th power-sum symmetric function.
In ghost coordinates,
$p_j\circ(w_1,w_2,\ldots)=(w_j,w_{2j},\ldots)$ in the presentation $W(\C) = \prod_{k\geq 1} \C$.  For an admissible $\Z$-set
$V$, let $C(V,W(\C))$ be the $\lambda$-ring of $W(\C)$-valued functions on
$V$ that are constant on $\Z$-orbits. For more details about $\lambda$-rings,
see \cite[Section~2]{Howe2024}.

Let $\varphi:V\to B$ be a morphism of admissible $\Z$-sets.
We recall the definitions of $\varphi^*$ and $\int_{V/B}$ below:
\begin{align*}
  (\varphi^*g)(v)
  &=p_{\deg(v)/\deg(\varphi(v))}\circ g(\varphi(v)),\\
  \left(\int_{V/B}f\right)(b)
  &=\sum_{|v|\in|V_b|} f(|v|)\bigl(t^{\deg_{V_b}(|v|)}\bigr),
\end{align*}
for $g\in C(B, W(\C))$ and $f\in C(V, W(\C))$.
For $w\in W(\C)$, the notation $w(t^d)$ denotes the substitution
$t\mapsto t^d$ in the presentation
$W(\C)=1+t\C[[t]]$.  Set $[V/B]:=\int_{V/B}1$.  If $V\to B$ admits a
section, then $[V/B]$ is invertible and the associated $W(\C)$-valued expectation
is
\[
  \E_{V/B}[f]:=\frac{\int_{V/B}f}{[V/B]}.
\]
In particular, for the final morphism $V\to\mathbf 1$, we write $\int_V$, $[V]$, and
$\E_V$.  The $k$th ghost coordinate of this expectation is an ordinary finite
average:
\begin{align}\label{eq:expectation-bridge}
    \gh_k\!\left(\E_V[f]\right)
  =
  \frac{1}{|V_k(1)|}
  \sum_{u\in V_k(1)}
  \gh_1\!\left(f_{V_k}(u)\right),
\end{align}
where $f_{V_k}$ is the canonical restriction of $f$ to $V_k$.

\subsection{Plethystic exponential and logarithm}

Let \(R\) be a pre-\(\lambda\)-ring.  We write
\[
  R[[\underline{t}_{\mathbf N}]]
  :=
  \varprojlim_N R[[t_1,\ldots,t_N]]
\]
for the ring of formal power series in countably many variables.  It is filtered
by total monomial degree: \(\Fil^d R[[\underline{t}_{\mathbf N}]]\) consists of series whose
monomials all have total degree at least \(d\).  Thus \(\Fil^1 R[[\underline{t}_{\mathbf N}]]\)
is the ideal of series with zero constant term, and
\(1+\Fil^1 R[[\underline{t}_{\mathbf N}]]\) is the set of series with constant term \(1\).

The plethystic \(\lambda\)-ring structure extends to \(R[[\underline{t}_{\mathbf N}]]\) by
requiring
\[
  p_i\circ t_j=t_j^i.
\]
Let \(h_k\) denote the \(k\)-th complete symmetric function.  The plethystic
exponential is the map
\[
  \Exp:\Fil^1 R[[\underline{t}_{\mathbf N}]]
  \longrightarrow
  1+\Fil^1 R[[\underline{t}_{\mathbf N}]],
  \qquad
  F\longmapsto \sum_{k\ge0} h_k\circ F.
\]
Here, the action $h_k\circ (-)$ can be obtained by expressing $h_k$ in terms
of the power sums $p_i$ using Newton's identities.
The filtration makes this sum well-defined coefficientwise.  Its inverse is denoted
\[
  \Log:
  1+\Fil^1 R[[\underline{t}_{\mathbf N}]]
  \longrightarrow
  \Fil^1 R[[\underline{t}_{\mathbf N}]].
\]
Thus, for \(H\in 1+\Fil^1 R[[\underline{t}_{\mathbf N}]]\) and \(Y\in R[[\underline{t}_{\mathbf N}]]\),
we use the usual plethystic power notation
\[
  H^Y:=\Exp\bigl(Y\,\Log H\bigr).
\]

\subsection{Motivic Euler products}

We recall the motivic Euler products below.
Let \(V\to B\) be a morphism of admissible \(\Z\)-sets.  Applying the above
with \(R=C(V,W(\C))\), for
\[
  H\in 1+\Fil^1 C(V,W(\C))[[\underline{t}_{\mathbf N}]]
\]
we define the motivic Euler product over \(V/B\) by
\[
  \prod_{V/B} H
  :=
  \Exp\left(
    \int_{V/B}\Log H
  \right)
  \in C(B,W(\C))[[\underline{t}_{\mathbf N}]].
\]
When \(B=\mathbf 1\), we write
\[
  \prod_V H
  :=
  \Exp\left(
    \int_V \Log H
  \right).
\]
The point-counting formula of \cite[Proposition~3.2.1]{BH} converts this motivic
product into a classical Euler product after a ghost coordinate is taken:
\[
  \gh_k\!\left(\prod_V H\right)
  =
  \prod_{|v|\in|V_k|}
  \gh_1\!\left(H(|v|)\right)
  \bigl(\underline{t}^{\deg(|v|)}\bigr).
\]

\subsection{Equidistribution}
We now recall the abstract equidistribution condition.

\begin{definition}[Bertucci--Howe equidistribution]\label{def:BH-equidist}
Let $A\to B$ be a morphism of admissible $\Z$-sets admitting a section, let
$I$ be a directed set, and for every $D\in I$ let $U_D$ be an admissible
$\Z$-set equipped with a morphism
$\ev_D:U_D\times B\to A$ over $B$.  For $k\ge1$ and a finite admissible
$\Z$-subset $B'\subseteq B_k$, define
\[
  \ev_{D,B'}:U_{D,k}(1)\longrightarrow
  \operatorname{Hom}_{B_k}(B',A_k),
  \qquad
  u\longmapsto\bigl(b\mapsto\ev_{D,k}(u,b)\bigr).
\]
The family $(U_D,\ev_D)_{D\in I}$ \emph{equidistributes on $A/B$} if, for
every $k\ge1$ and every nonempty finite admissible $B'\subseteq B_k$,
\[
  \lim_{D\in I}(\ev_{D,B'})_*\mu_{U_{D,k}(1)}
  =
  \mu_{\operatorname{Hom}_{B_k}(B',A_k)},
\]
where $\mu_S$ denotes the uniform probability measure on a finite set $S$.
\end{definition}

\subsection{Asymptotic \texorpdfstring{\(\sigma\)}{sigma}-moment generating functions}
The following joint form is Theorem~4.2.1 and Remark~4.2.2 of
\cite{BH}.
\begin{theorem}[Bertucci--Howe, joint form]\label{thm:BH-joint}
Let $I$ be a directed set and assume that $(U_D,\ev_D)_{D\in I}$
equidistributes on $A/B$ in the sense of
Definition~\ref{def:BH-equidist}.  Let
\[
  \chi^{(1)},\ldots,\chi^{(n)}\in C(A,W(\C))
\]
be local statistics, and set
\[
  X_D^{(r)}:=\int_{U_D\times B/U_D}\ev_D^*\chi^{(r)}.
\]
Then, for
\[
  H_r\in 1+\Fil^1\Z[[t_{\mathbf N}^{(r)}]]
  \qquad (1\le r\le n),
\]
where the variables $t_{\mathbf N}^{(r)}$ are independent for different $r$,
\[
  \lim_{D\in I}
  \E_{U_D}\left[
    \prod_{r=1}^n H_r^{X_D^{(r)}}
  \right]
  =
  \prod_B
  \E_{A/B}\left[
    \prod_{r=1}^n H_r^{\chi^{(r)}}
  \right].
\]
\end{theorem}

\subsection{Specialization to the \texorpdfstring{\(\kappa\)}{kappa}-power-free
superelliptic family}
\label{subsec:superelliptic-specialization}

We now take \(I=\N\), with \(D\) the polynomial degree.  Let
\(B=\PP^1(\overline{\F}_q)\), equipped with the geometric Frobenius action.  We
write \(\overline U_D\) for the admissible \(\Z\)-set of degree-\(D\),
\(\kappa\)-power-free polynomials with coefficients in \(\overline{\F}_q\).
Its fixed-point set \(\overline U_D(1)\) is the set \(U_D\) introduced in
Section~\ref{sec:introduction}, and more generally
\[
  \overline U_{D,k}(1)
  =
  \left\{
    F\in\F_{q^k}[x]:
    \deg F=D
    \text{ and }F\text{ is \(\kappa\)-power-free}
  \right\}.
\]

Define an admissible \(\Z\)-set \(A\) over \(B\) by prescribing its geometric
fibers.  For a finite geometric point
\(P\in\A^1(\overline{\F}_q)\), write
\[
  \cO_P:=\cO_{\PP^1_{\overline{\F}_q},P},
  \qquad
  \frakm_P\subset\cO_P
\]
for its local ring and maximal ideal, and put
\[
  A_P
  :=
  \left(\cO_P/\frakm_P^\kappa\right)\setminus\{0\},
  \qquad
  \ev_D(F,P):=F\bmod\frakm_P^\kappa.
\]
At \(\infty=[1:0]\), put
\[
  A_\infty=\overline{\F}_q^\times,
  \qquad
  \ev_D(F,\infty):=\operatorname{lc}(F).
\]
The constant jets \(1\in A_P\) at finite \(P\), together with
\(1\in A_\infty\), define a Frobenius-equivariant section \(B\to A\), as
required by Theorem~\ref{thm:BH-joint}.
Indeed, after homogenizing \(F\) and trivializing by \(x_0^D\) near infinity,
the resulting value at infinity is the leading coefficient of \(F\).  Thus the
condition that \(F\) be \(\kappa\)-power-free is exactly the condition that every
finite local \(\kappa\)-jet lie in the nonzero geometric fiber \(A_P\), while the
condition at infinity records \(\deg F=D\).

We next introduce the corresponding fibers over Frobenius orbits.  For a
geometric point \(P\in B\), let
\[
  |P|
  :=
  \left\{
    P,\operatorname{Frob}_q(P),
    \operatorname{Frob}_q^2(P),\ldots
  \right\}
\]
be its Frobenius orbit, and put
\[
  e_P:=\deg(P):=\#|P|.
\]
The orbit \(|P|\) is identified with the corresponding closed point of
\(\PP^1_{\F_q}\).  Define
\[
  A_{|P|}
  :=
  \operatorname{Hom}_B(|P|,A).
\]
This is the orbitwise fiber in the equidistribution space of
Definition~\ref{def:BH-equidist}.  See \cite[\S4.1]{BH}.

An element of \(A_{|P|}\) is a Frobenius-equivariant section of \(A\to B\) over
the orbit \(|P|\).  Evaluation at the chosen representative \(P\) identifies
\[
  A_{|P|}
  \cong
  A_P^{\operatorname{Frob}_q^{e_P}}.
\]
If \(P\) is finite, this gives a canonical identification
\[
  A_{|P|}
  \cong
  \left(
    \cO_{\A^1_{\F_q},|P|}/\frakm_{|P|}^{\kappa}
  \right)\setminus\{0\}.
\]
After choosing a uniformizer at the closed point \(|P|\), one obtains a
noncanonical identification
\[
  A_{|P|}
  \cong
  \left(
    \F_{q^{e_P}}[\varepsilon]/(\varepsilon^\kappa)
  \right)\setminus\{0\}.
\]
This is the closed-point local space denoted \(A_{|P|}\) in
\cite[Example~4.1.2 and Proposition~5.1.1]{BH}.  More generally, for an orbit
in \(B_k\), we write
\[
  A_{k,|P|}
  :=
  \operatorname{Hom}_{B_k}(|P|,A_k).
\]

For a point \((P,\alpha)\in A\), let
\[
  d(P,\alpha):=\deg_A(P,\alpha)
\]
be the size of its Frobenius orbit in the admissible \(\Z\)-set \(A\).  Since
\(A\to B\) is Frobenius-equivariant, one has \(e_P\mid d(P,\alpha)\).  If
\(P\) is finite, let \(\overline\alpha\) be the image of \(\alpha\) under the
canonical reduction map
\[
  \cO_P/\frakm_P^\kappa
  \longrightarrow
  \cO_P/\frakm_P
  \cong
  \overline{\F}_q.
\]
Thus \(\overline\alpha\) is the constant term, or residue, of the jet.  When
\(\alpha=\ev_D(F,P)\), it is simply \(F(P)\).

For each \(r\ge1\), define a function
\[
  \chi^{(r)}:A\longrightarrow\Z_{\ge0}
\]
by
\[
  \chi^{(r)}(P,\alpha)
  :=
  \begin{cases}
    \displaystyle
    e_P\,\#\{y\in\F_{q^r}:y^m=\overline\alpha\},
    & P\ne\infty,\quad d(P,\alpha)=e_P,\quad e_P\mid r,\\
    0,&\text{otherwise}.
  \end{cases}
\]
In particular, \(\chi^{(r)}(\infty,\alpha)=0\).
When the first case occurs, the condition
\(d(P,\alpha)=e_P\) implies
\(\overline\alpha\in\F_{q^{e_P}}\subseteq\F_{q^r}\), so the displayed
solution set is well-defined.  The function \(\chi^{(r)}\) is constant on
Frobenius orbits.  Indeed, Frobenius preserves \(e_P\) and
\(d(P,\alpha)\), and its action on \(\F_{q^r}\) carries the equation
\(y^m=\overline\alpha\) bijectively to the corresponding conjugate equation.
Here \([1]\in W(\C)\) denotes the one-point class, and \(\nu[1]\) denotes its
\(\nu\)-fold sum.
Consequently,
\[
  \chi_W^{(r)}:A\longrightarrow W(\C),
  \qquad
  (P,\alpha)\longmapsto\chi^{(r)}(P,\alpha)[1],
\]
is a well-defined element of \(C(A,W(\C))\), as required by
Theorem~\ref{thm:BH-joint}.  The values on points with
\(d(P,\alpha)>e_P\) have been set to zero.  These points do not occur in the
first-ghost fiber calculations below.

For an orbitwise section
\[
  \alpha\in A_{|P|}=\operatorname{Hom}_B(|P|,A),
\]
define
\[
  \chi^{(r)}(|P|,\alpha)
  :=
  \chi^{(r)}\bigl(P,\alpha(P)\bigr),
\]
where \(P\) is any geometric representative of \(|P|\).  This is independent
of the representative by Frobenius invariance.  Equivalently,
\[
  \chi^{(r)}(|P|,\alpha)
  =
  \sum_{Q\in|P|\cap\A^1(\F_{q^r})}
  \#\{y\in\F_{q^r}:y^m=\overline{\alpha(Q)}\}.
\]
In particular, set \(\chi^{(r)}(|\infty|,-)=0\).  The explicit value of this
orbitwise statistic will be computed in Section~\ref{sec:local-factors}.

For \(F\in U_D\), Frobenius equivariance gives an orbitwise evaluation
\[
  \ev_D(F,|P|)
  :=
  \bigl(Q\longmapsto\ev_D(F,Q)\bigr)
  \in A_{|P|}.
\]
Because \(F\) is fixed by Frobenius,
\[
  d\bigl(P,\ev_D(F,P)\bigr)=e_P.
\]
Hence \(\chi^{(r)}(|P|,\ev_D(F,|P|))\) is precisely the total contribution
of the orbit \(|P|\) to the affine point count over \(\F_{q^r}\).

The global statistic appearing in Theorem~\ref{thm:BH-joint} is
\begin{equation}\label{eq:Xd-superelliptic}
  X_D^{(r)}
  :=
  \int_{\overline U_D\times B/\overline U_D}
  \ev_D^*\chi_W^{(r)}.
\end{equation}
For \(F\in U_D=\overline U_D(1)\), the point \((F,P)\) and its image
\(\ev_D(F,P)\) have the same orbit degree \(e_P\).  Thus the Adams operation
in the pullback \(\ev_D^*\) is \(p_1\), and integration over fibers gives
\[
  X_D^{(r)}(F)
  =
  \sum_{|P|\in|\PP^1_{\F_q}|}
  \left(
    \chi^{(r)}\bigl(|P|,\ev_D(F,|P|)\bigr)[1]
  \right)(t^{e_P}).
\]
The affine point count is obtained by summing the same orbit multiplicities:
\[
  N_r(F)
  =
  \sum_{|P|\in|\PP^1_{\F_q}|}
  \chi^{(r)}\bigl(|P|,\ev_D(F,|P|)\bigr).
\]
Indeed, a finite orbit \(|P|\) contributes only when \(e_P\mid r\).  In that
case the factor \(e_P\) in the definition of \(\chi^{(r)}\) accounts for all
geometric points in the orbit.  The orbit at infinity contributes zero.

Fix \(n\), and choose in Theorem~\ref{thm:BH-joint}
\[
  H_r=1+h_1^{(r)}
  \qquad (1\le r\le n).
\]
Once Proposition~\ref{prop:equidist} supplies equidistribution for this family,
Theorem~\ref{thm:BH-joint} gives the specialization \(H=1+h_1\) from
\cite[Remark~4.2.2]{BH}, applied jointly to the \(n\) independent alphabets,
whose first complete symmetric functions are \(h_1^{(r)}\):
\begin{equation}\label{eq:BH-plugged-superelliptic}
  \lim_{D\to\infty}
  \E_{\overline U_D}
  \left[
    \prod_{r=1}^n
    \bigl(1+h_1^{(r)}\bigr)^{X_D^{(r)}}
  \right]
  =
  \prod_B
  \E_{A/B}
  \left[
    \prod_{r=1}^n
    \bigl(1+h_1^{(r)}\bigr)^{\chi_W^{(r)}}
  \right].
\end{equation}

We now prove the required equidistribution statement.

\begin{proposition}\label{prop:equidist}
The evaluation maps
\[
  \ev_D:
  \overline U_D\times\PP^1(\overline{\F}_q)
  \longrightarrow A
\]
equidistribute on \(A/\PP^1(\overline{\F}_q)\) as \(D\to\infty\).
\end{proposition}

\begin{proof}
We apply \cite[Proposition~5.1.1]{BH} to the admissible \(\Z\)-set of
homogeneous binary forms
\[
  S_D(\overline{\F}_q)
  =
  H^0\left(\PP^1_{\overline{\F}_q},\cO(D)\right),
\]
using the geometric subsets \(A_P\) defined above.  For
\(P\in\A^1(\overline{\F}_q)\), take \(j_P=1\) and
\(M_P=\kappa-1\).  At \(\infty=[1:0]\), take \(j_\infty=0\) and
\(M_\infty=0\).

For the Taylor-condition hypothesis in
\cite[Proposition~5.1.1]{BH}, take
\(Y=\A^1_{\F_q}\subset\PP^1_{\F_q}\) and the identity quotient of
\(\iota^*\Omega^1_{\PP^1_{\F_q}}\).  At a finite geometric point \(P\), the
resulting condition says that the image of \(G/x_1^D\) in
\[
  \cO_{\PP^1_{\overline{\F}_q},P}/\frakm_P^2
\]
is nonzero.  Since \(\kappa\ge2\), this implies that its image in
\[
  \cO_{\PP^1_{\overline{\F}_q},P}/\frakm_P^\kappa
\]
is nonzero, and hence lies in \(A_P\).  The point \(\infty\) is a finite
exceptional set, which is permitted in
\cite[Proposition~5.1.1]{BH}.  Therefore the corresponding homogeneous
evaluation maps equidistribute.

It remains to identify this homogeneous family with the affine family above.
For
\[
  F(x)=a_Dx^D+\cdots+a_0,
  \qquad a_D\ne0,
\]
put
\[
  \widetilde F(x_0,x_1)
  :=
  x_1^D F(x_0/x_1).
\]
This gives a Frobenius-equivariant bijection between \(\overline U_D\) and the
corresponding homogeneous family.  If
\(P\in\A^1(\overline{\F}_q)\), then \(x_1\) is a unit at \(P\), and
\[
  \frac{\widetilde F}{x_1^D}=F(x_0/x_1).
\]
Hence the image of \(\widetilde F/x_1^D\) in
\(\cO_{\PP^1_{\overline{\F}_q},P}/\frakm_P^\kappa\) is exactly the image of
\(F\) modulo \(\frakm_P^\kappa\).  At infinity, the image of
\(\widetilde F/x_0^D\) in
\[
  \cO_{\PP^1_{\overline{\F}_q},\infty}/\frakm_\infty
\]
is the leading coefficient \(a_D\).

Thus the homogeneous evaluation map is identified with \(\ev_D\).  The same
Frobenius-equivariant bijection identifies the \(k\)-fixed-point sets and their
uniform probability measures.  Therefore the equidistribution statement of
\cite[Proposition~5.1.1]{BH} transfers to the evaluation maps above.
\end{proof}

\section{Local computation for superelliptic curves}
\label{sec:local-factors}

We retain the notation of
Subsection~\ref{subsec:superelliptic-specialization}.  Fix a finite geometric
point \(P\in\A^1(\overline{\F}_q)\), with Frobenius orbit \(|P|\), degree
\(e_P\), and orbit fiber \(A_{|P|}\).  The local statistic on the orbit
fiber is induced by the function defined on all of \(A\):
\[
  \chi^{(r)}(|P|,\alpha)
  =
  \chi^{(r)}\bigl(P,\alpha(P)\bigr),
  \qquad \alpha\in A_{|P|}.
\]
Using the local description given there, write
\[
  \alpha
  =
  u_0+u_1\varepsilon+\cdots+u_{\kappa-1}\varepsilon^{\kappa-1},
  \qquad
  u_i\in\F_{q^{e_P}},
\]
not all zero.  Let \(v_{|P|}(\alpha)\) be the least \(i\) for which
\(u_i\ne0\).

\begin{lemma}\label{lem:local-contribution}
If \(e_P\nmid r\), then \(\chi^{(r)}(|P|,\alpha)=0\).  If \(e_P\mid r\),
then
\[
  \chi^{(r)}(|P|,\alpha)
  =
  \begin{cases}
    e_P,
    & v_{|P|}(\alpha)\ge1,\\[4pt]
    e_Pg_r\,\mathbf 1_{\{u_0\in H_r^{(e_P)}\}},
    & v_{|P|}(\alpha)=0,
  \end{cases}
\]
where
\[
  H_r^{(e_P)}
  :=
  (\F_{q^r}^{\times})^m\cap\F_{q^{e_P}}^{\times}.
\]
\end{lemma}

\begin{proof}
The orbit \(|P|\) consists of \(e_P\) geometric points and is contained in
\(\A^1(\F_{q^r})\) exactly when \(e_P\mid r\).  Thus, if \(e_P\nmid r\), it
contributes no \(\F_{q^r}\)-points of the base.

Assume \(e_P\mid r\).  If \(v_{|P|}(\alpha)\ge1\), then the residue of
\(\alpha\) vanishes at every point of \(|P|\).  Above each of these \(e_P\)
points, the equation \(y^m=0\) has exactly one solution, so the total
contribution is \(e_P\).

If \(v_{|P|}(\alpha)=0\), then the residue at \(P\) is
\(u_0\in\F_{q^{e_P}}^{\times}\), while the residues at the other points of the
orbit are its Frobenius conjugates.  The equation \(y^m=u_0\) has
\(g_r=\gcd(m,q^r-1)\) solutions in \(\F_{q^r}\) precisely when
\(u_0\in(\F_{q^r}^{\times})^m\).  Since
\(u_0\in\F_{q^{e_P}}^{\times}\), this is equivalent to
\(u_0\in H_r^{(e_P)}\).  The subgroup of \(m\)-th powers is Frobenius-stable,
so the same condition holds at every point of the orbit.  Hence the total
contribution is \(e_Pg_r\) if \(u_0\in H_r^{(e_P)}\), and zero otherwise.
\end{proof}

The following calculation in the cyclic group \(\F_{q^{e_P}}^{\times}\)
controls the dependence between different extension degrees.

\begin{lemma}\label{lem:index}
For \(e_P\mid r\), the subgroup \(H_r^{(e_P)}\) has index
\[
  [\F_{q^{e_P}}^{\times}:H_r^{(e_P)}]
  =
  \delta_r^{(e_P)}
  =
  \frac{g_r}{\gcd(g_r,c_r^{(e_P)})}.
\]
Equivalently,
\[
  \frac{|H_r^{(e_P)}|}{q^{e_P}-1}
  =
  \frac1{\delta_r^{(e_P)}}
  =
  \frac{\gamma_r^{(e_P)}}{g_r}.
\]
\end{lemma}

\begin{proof}
Choose a generator \(\zeta\) of \(\F_{q^r}^{\times}\), and put
\[
  c=c_r^{(e_P)}=\frac{q^r-1}{q^{e_P}-1}.
\]
Then
\[
  \F_{q^{e_P}}^{\times}=\langle\zeta^c\rangle,
  \qquad
  (\F_{q^r}^{\times})^m=\langle\zeta^{g_r}\rangle.
\]
An element \(\zeta^{cj}\) belongs to \(\langle\zeta^{g_r}\rangle\) if and only
if \(g_r\mid cj\), or equivalently
\[
  \frac{g_r}{\gcd(g_r,c)}\mid j.
\]
This gives the stated index.
\end{proof}

\begin{proposition}\label{prop:pattern-counts}
For \(S\subseteq R_{e_P}\), the number of units in
\(\F_{q^{e_P}}^{\times}\) whose membership pattern is exactly \(S\) is
\[
  C_S^{(e_P)}
  =
  \sum_{S\subseteq T\subseteq R_{e_P}}
  (-1)^{|T|-|S|}
  \frac{q^{e_P}-1}{\lcmop_{r\in T}\delta_r^{(e_P)}},
  \qquad
  \lcmop_{\varnothing}:=1.
\]
\end{proposition}

\begin{proof}
For \(T\subseteq R_{e_P}\), put
\[
  H_T:=\bigcap_{r\in T}H_r^{(e_P)},
  \qquad
  H_{\varnothing}:=\F_{q^{e_P}}^{\times}.
\]
A unit has membership pattern exactly \(S\) if it lies in \(H_S\) and avoids
\(H_r^{(e_P)}\) for every \(r\in R_{e_P}\setminus S\).  Inclusion--exclusion
therefore gives
\[
  C_S^{(e_P)}
  =
  \sum_{S\subseteq T\subseteq R_{e_P}}
  (-1)^{|T|-|S|}|H_T|.
\]
Since \(\F_{q^{e_P}}^{\times}\) is cyclic, the intersection of its subgroups
of indices \(\delta_r^{(e_P)}\), for \(r\in T\), is the subgroup of index
\(\lcmop_{r\in T}\delta_r^{(e_P)}\).  Hence
\[
  |H_T|
  =
  \frac{q^{e_P}-1}{\lcmop_{r\in T}\delta_r^{(e_P)}}.
\]
Substitution proves the result.
\end{proof}

\begin{theorem}\label{thm:local-factor}
The local probability generating function attached to \(|P|\) is
\[
  \begin{aligned}
    L_{e_P}(\mathbf w)
    &:={}
    \E_{\alpha\in A_{|P|}}
    \prod_{r=1}^n w_r^{\chi^{(r)}(|P|,\alpha)} \\
    &={}
    \alpha_{e_P}\prod_{r\in R_{e_P}}w_r^{e_P}
    +(1-\alpha_{e_P})
    \sum_{S\subseteq R_{e_P}}
    \frac{C_S^{(e_P)}}{q^{e_P}-1}
    \prod_{r\in S}w_r^{e_Pg_r}.
  \end{aligned}
\]
In particular, the local factor depends only on the orbit degree \(e_P\).
\end{theorem}

\begin{proof}
Let \(\alpha\) be uniformly distributed on \(A_{|P|}\).  The positive-valuation
classes are exactly the nonzero classes with zero constant coefficient.  There
are \(q^{e_P(\kappa-1)}-1\) such classes, while
\(|A_{|P|}|=q^{e_P\kappa}-1\).  Consequently,
\[
  \Prb\bigl(v_{|P|}(\alpha)\ge1\bigr)
  =
  \frac{q^{e_P(\kappa-1)}-1}{q^{e_P\kappa}-1}
  =
  \alpha_{e_P}.
\]
If \(v_{|P|}(\alpha)\ge1\), Lemma~\ref{lem:local-contribution} gives
\(\chi^{(r)}(|P|,\alpha)=e_P\) for \(r\in R_{e_P}\), and zero otherwise.
This case contributes
\[
  \alpha_{e_P}\prod_{r\in R_{e_P}}w_r^{e_P}.
\]

It remains to consider \(v_{|P|}(\alpha)=0\).  For every fixed
\(u\in\F_{q^{e_P}}^{\times}\), there are exactly
\(q^{e_P(\kappa-1)}\) choices of the remaining coefficients of \(\alpha\).
Thus, conditional on \(v_{|P|}(\alpha)=0\), the constant coefficient \(u_0\)
is uniform on \(\F_{q^{e_P}}^{\times}\).  Lemma~\ref{lem:local-contribution}
therefore gives
\[
  \begin{aligned}
    &\E\left[
      \prod_{r=1}^n w_r^{\chi^{(r)}(|P|,\alpha)}
      \,\middle|\,
      v_{|P|}(\alpha)=0
    \right] \\
    &\qquad={}
    \frac1{q^{e_P}-1}
    \sum_{u\in\F_{q^{e_P}}^{\times}}
    \prod_{r\in R_{e_P}}
    w_r^{e_Pg_r\mathbf 1_{\{u\in H_r^{(e_P)}\}}} \\
    &\qquad={}
    \sum_{S\subseteq R_{e_P}}
    \frac{C_S^{(e_P)}}{q^{e_P}-1}
    \prod_{r\in S}w_r^{e_Pg_r}.
  \end{aligned}
\]
Combining the two cases proves the formula.
\end{proof}

\section{The joint distribution}
\label{sec:main-proof}

We apply Theorem~\ref{thm:BH-joint} to the Witt-valued statistics
\(\chi_W^{(r)}\) constructed in
Subsection~\ref{subsec:superelliptic-specialization}. Taking the first ghost
coordinate and specializing the symmetric-function variables gives an
ordinary probability generating function.

\begin{lemma}\label{lem:ordinary-specialization}
For fixed \(n\),
\[
  \lim_{D\to\infty}
  \E_{F\in U_D}\prod_{r=1}^n w_r^{N_r(F)}
  =
  \prod_{|P|\in S_n}
  \E_{\alpha\in A_{|P|}}
  \prod_{r=1}^n w_r^{\chi^{(r)}(|P|,\alpha)},
\]
where
\[
  S_n:=\{\,|P|\in|\A^1_{\F_q}|:e_P\le n\,\}.
\]
\end{lemma}

\begin{proof}
We start from the identity \eqref{eq:BH-plugged-superelliptic}.  Let
\(p_i^{(r)}\) denote the \(i\)-th power sum in the \(r\)-th independent
alphabet.  On the completed symmetric-function ring in these \(n\) alphabets,
define
\[
  \Theta_n(p_i^{(r)})
  =
  \begin{cases}
    w_r-1,&1\le i\le n,\\
    0,&i>n.
  \end{cases}
\]
The source is completed by total symmetric degree and the target is
\(\C[[w_1-1,\ldots,w_n-1]]\).  If \(\Fil^d\) denotes the source filtration and
\(\mathfrak n=(w_1-1,\ldots,w_n-1)\), then
\[
  \Theta_n(\Fil^d)\subseteq\mathfrak n^{\lceil d/n\rceil}.
\]
Hence \(\Theta_n\) is continuous.  We may therefore apply
\(\Theta_n\circ\gh_1\) coefficientwise to the completed series and to their
limits.

We first record the local specialization.  Fix \(e,r\ge1\) and
\(\nu\in\Z_{\ge0}\), and put
\[
  \ell_r:=\Log(1+h_1^{(r)}).
\]
By the definition of the pre-\(\lambda\) power,
\[
  (1+h_1^{(r)})^{\nu[1](t^e)}
  =
  \Exp\bigl(\nu[1](t^e)\ell_r\bigr).
\]
By \cite[Lemma~2.5.4]{Howe2024},
\[
  \Exp(x)=\prod_{i\ge1}\exp\left(\frac{p_i\circ x}{i}\right).
\]
Taking the first ghost coordinate and using multiplicativity of the Adams
operations gives
\[
  \begin{aligned}
    \gh_1\Exp\bigl(\nu[1](t^e)\ell_r\bigr)
    &={}
    \prod_{i\ge1}
    \exp\left(
      \frac{\nu\,\gh_1(p_i\circ[1](t^e))\,
      \gh_1(p_i\circ \ell_r)}{i}
    \right).
  \end{aligned}
\]
By the substitution formula for Witt ghost coordinates,
\[
  \gh_1(p_i\circ[1](t^e))
  =
  \gh_i([1](t^e))
  =
  \begin{cases}
    e,&e\mid i,\\
    0,&e\nmid i.
  \end{cases}
\]
Thus only \(i=ej\) contribute, and the preceding product is
\[
  \prod_{j\ge1}
  \exp\left(\nu\frac{\gh_1(p_{ej}\circ \ell_r)}{j}\right)
  =
  \gh_1\Exp\bigl(\nu\,p_e\circ \ell_r\bigr).
\]
Since Adams operations commute with \(\Log\),
\[
  p_e\circ \ell_r=\Log(1+p_e^{(r)}).
\]
Consequently,
\[
  \gh_1\left((1+h_1^{(r)})^{\nu[1](t^e)}\right)
  =(1+p_e^{(r)})^\nu,
\]
and therefore
\[
  (\Theta_n\circ\gh_1)\left((1+h_1^{(r)})^{\nu[1](t^e)}\right)
  =
  \begin{cases}
    w_r^\nu,&e\le n,\\
    1,&e>n.
  \end{cases}
\]

We apply this calculation first to the left-hand side of
\eqref{eq:BH-plugged-superelliptic}.  For \(F\in U_D\),
\[
  X_D^{(r)}(F)
  =
  \sum_{|P|\in|\PP^1_{\F_q}|}
  \left(
    \chi^{(r)}(|P|,\ev_D(F,|P|))[1]
  \right)(t^{e_P}).
\]
Additivity of the pre-\(\lambda\) exponent and the preceding local calculation
give
\[
  \begin{aligned}
    &(\Theta_n\circ\gh_1)
    \left[
      \prod_{r=1}^n(1+h_1^{(r)})^{X_D^{(r)}(F)}
    \right] \\
    &\qquad={}
    \prod_{|P|\in|\PP^1_{\F_q}|}
    \prod_{r=1}^n
    w_r^{\chi^{(r)}(|P|,\ev_D(F,|P|))}
    =
    \prod_{r=1}^n w_r^{N_r(F)}.
  \end{aligned}
\]
By the \(k=1\) case of \cite[Lemma~2.4.4]{BH}, the first ghost coordinate of
the expectation over \(\overline U_D\) is the ordinary uniform average over
\(\overline U_D(1)=U_D\).  Thus the specialized left-hand side is
\[
  \lim_{D\to\infty}
  \E_{F\in U_D}\prod_{r=1}^n w_r^{N_r(F)}.
\]

It remains to specialize the right-hand side.  Put
\[
  \mathcal R
  :=
  \prod_B
  \E_{A/B}
  \left[
    \prod_{r=1}^n(1+h_1^{(r)})^{\chi_W^{(r)}}
  \right].
\]
By \cite[Proposition~3.2.1]{BH},
\[
  \begin{aligned}
    \gh_1(\mathcal R)
    =
    \prod_{|P|\in|B|}
    \gh_1\left(
      \left.
      \E_{A/B}
      \left[
        \prod_{r=1}^n
        (1+p_{e_P}^{(r)})^{\chi_W^{(r)}}
      \right]
      \right|_P
    \right),
  \end{aligned}
\]
where \(P\) is any geometric representative of \(|P|\).  By
\cite[Lemma~2.4.7]{BH}, evaluation at \(P\) gives the Howe expectation on the
geometric fiber \(A_P\), equipped with its fiber action.  Its fixed-point set is
naturally
\[
  A_P(1)\cong\operatorname{Hom}_B(|P|,A)=A_{|P|}.
\]
Applying the \(k=1\) case of \cite[Lemma~2.4.4]{BH}, and using the explicit
construction of \(\chi_W^{(r)}\), therefore gives
\[
  \begin{aligned}
    \gh_1(\mathcal R)
    =
    \prod_{|P|\in|B|}
    \E_{\alpha\in A_{|P|}}
    \left[
      \gh_1\left(
        \prod_{r=1}^n
        (1+p_{e_P}^{(r)})^{\chi^{(r)}(|P|,\alpha)[1]}
      \right)
    \right].
  \end{aligned}
\]
For every \(\nu\ge0\),
\[
  \gh_1\left((1+p_{e_P}^{(r)})^{\nu[1]}\right)
  =(1+p_{e_P}^{(r)})^\nu.
\]
Hence, after applying \(\Theta_n\), the orbit at infinity and every orbit with
\(e_P>n\) contribute \(1\), and
\[
  (\Theta_n\circ\gh_1)(\mathcal R)
  =
  \prod_{|P|\in S_n}
  \E_{\alpha\in A_{|P|}}
  \prod_{r=1}^n w_r^{\chi^{(r)}(|P|,\alpha)}.
\]
Applying the continuous map \(\Theta_n\circ\gh_1\) to
\eqref{eq:BH-plugged-superelliptic} gives the asserted equality in the completed
target.  Both sides are polynomials in \(\mathbf w\), so this is also equality
as ordinary polynomials and proves the lemma.
\end{proof}

\begin{proof}[Proof of Theorem~\ref{thm:main-pgf}]
By Proposition~\ref{prop:equidist}, the Bertucci--Howe theorem applies to the
superelliptic family.  Lemma~\ref{lem:ordinary-specialization} and
Theorem~\ref{thm:local-factor} show that every finite orbit of degree \(e\)
contributes the same factor \(L_e(\mathbf w)\).  There are \(a_e\) such closed
points, and only \(e\le n\) can contribute to the first \(n\) coordinates.
Therefore
\[
  \lim_{D\to\infty}
  \E_{F\in U_D}\prod_{r=1}^n w_r^{N_r(F)}
  =
  \prod_{e=1}^nL_e(\mathbf w)^{a_e}.
\]
For fixed \(q,m,n\),
\[
  (N_1(F),\ldots,N_n(F))
  \in
  \prod_{r=1}^n\{0,\ldots,g_rq^r\}.
\]
Thus convergence of the pgfs is coefficientwise convergence of the probability
masses, which is equivalent to convergence in distribution on this finite state
space.
\end{proof}

\begin{proof}[Proof of Corollary~\ref{cor:blocks}]
For one degree-\(e\) block, the construction of \(X_e\) gives
\[
  \E\prod_{r=1}^n w_r^{X_e^{(r)}}
  =\alpha_e\prod_{r\in R_e}w_r^e
  +\frac{1-\alpha_e}{q^e-1}
    \sum_{u\in\F_{q^e}^{\times}}
    \prod_{r\in R_e}
    w_r^{eg_r\mathbf 1_{\{u\in H_r^{(e)}\}}}.
\]
Grouping the units by membership pattern shows that this pgf is
\(L_e(\mathbf w)\).  Hence the pgf of the independent block sum
\(\sum_{e=1}^n\sum_{j=1}^{a_e}X_{e,j}\) is
\[
  \prod_{e=1}^nL_e(\mathbf w)^{a_e}=G_n(\mathbf w).
\]
By Theorem~\ref{thm:main-pgf}, this is also the pgf of \(\mathbf N^*\).
A multivariate pgf uniquely determines the law of a finitely supported
\(\Z_{\ge0}^n\)-valued random vector, proving the assertion.
\end{proof}

\section{Moments and covariances}\label{sec:moments}

We now extract moment information from the independent closed-point block
description.
\subsection{Expectation and covariance}

For \(e\mid\gcd(r,s)\), put
\[
  \sigma_{r,s}^{(e)}
  :=e^2(1-\alpha_e)
  \left[
    \begin{aligned}
      &g_r g_s
      \left(
        \frac1{\lcmop(\delta_r^{(e)},\delta_s^{(e)})}
        -\frac1{\delta_r^{(e)}\delta_s^{(e)}}
      \right)\\[4pt]
      &{}+\alpha_e(1-\gamma_r^{(e)})(1-\gamma_s^{(e)})
    \end{aligned}
  \right].
\]

\Needspace{6\baselineskip}
\begin{theorem}[Means and covariances]
\label{thm:mean-covariance}
For \(r,s\ge1\),
\[
  \lim_{D\to\infty}
  \E_{F\in U_D}[N_r(F)]
  =
  \sum_{e\mid r}a_e e\left(\alpha_e+(1-\alpha_e)\gamma_r^{(e)}\right).
\]
Equivalently,
\[
  \lim_{D\to\infty}
  \E_{F\in U_D}[N_r(F)]
  =
  q^r+
  \sum_{e\mid r}
  a_e e\,(1-\alpha_e)\bigl(\gamma_r^{(e)}-1\bigr).
\]
Moreover,
\[
  \lim_{D\to\infty}
  \Cov_{F\in U_D}(N_r(F),N_s(F))
  =
  \sum_{e\mid\gcd(r,s)}a_e\sigma_{r,s}^{(e)}.
\]
\end{theorem}

\begin{proof}
Choose \(n\ge\max\{r,s\}\) in Corollary~\ref{cor:blocks}.  For the expectation,
fix \(e\mid r\), suppress the block index \(j\), and write
\(I_r=\mathbf 1_{\{U_e\in H_r^{(e)}\}}\).  Then
\[
  X_e^{(r)}
  =
  e\eta_e+e g_r(1-\eta_e)I_r.
\]
Since \(\eta_e\) and \(U_e\) are independent,
\(\E[\eta_e]=\alpha_e\), \(\E[1-\eta_e]=(1-\alpha_e)\), and
\[
  \E[I_r]=\frac{|H_r^{(e)}|}{q^e-1}=\frac1{\delta_r^{(e)}}.
\]
Therefore
\[
  \E[X_e^{(r)}]
  =
  e\alpha_e+
  e g_r(1-\alpha_e)\frac1{\delta_r^{(e)}}
  =
  e\left(\alpha_e+(1-\alpha_e)\gamma_r^{(e)}\right).
\]
Summing over the \(a_e\) degree-\(e\) blocks gives the first formula.  The second
formula follows from \(\sum_{e\mid r}ea_e=q^r\).

For covariance, only degree-\(e\) blocks with \(e\mid r\) and \(e\mid s\) can
contribute to both coordinates.  Fix such an \(e\), and write
\(I_s=\mathbf 1_{\{U_e\in H_s^{(e)}\}}\).  Since
\(\eta_e(1-\eta_e)=0\),
\[
  X_e^{(r)}X_e^{(s)}
  =
  e^2\eta_e+
  e^2 g_r g_s(1-\eta_e)I_r I_s.
\]
The group \(\F_{q^e}^{\times}\) is cyclic, and the subgroups
\(H_r^{(e)}\) and \(H_s^{(e)}\) have indices
\(\delta_r^{(e)}\) and \(\delta_s^{(e)}\).  Hence their intersection has index
\(\lcmop(\delta_r^{(e)},\delta_s^{(e)})\), so
\[
  \E[I_r I_s]
  =
  \frac1{\lcmop(\delta_r^{(e)},\delta_s^{(e)})}.
\]
Thus
\[
  \E[X_e^{(r)}X_e^{(s)}]
  =
  e^2
  \left(
    \alpha_e+
    (1-\alpha_e)
    \frac{g_r g_s}
    {\lcmop(\delta_r^{(e)},\delta_s^{(e)})}
  \right).
\]
Subtracting
\(\E[X_e^{(r)}]\E[X_e^{(s)}]\) gives \(\sigma_{r,s}^{(e)}\).  Distinct
closed-point blocks are independent, so all cross-block covariances vanish.
Therefore
\[
  \Cov(N_r^*,N_s^*)
  =
  \sum_{e\mid\gcd(r,s)}a_e\sigma_{r,s}^{(e)}.
\]
Moment convergence then gives the desired limiting covariance.
\end{proof}

\subsection{Computation of moments}

Fix \(r\ge1\) and choose any \(n\ge r\) in Corollary~\ref{cor:blocks}.  The
limiting variable \(N_r^*\) is the sum of the independent coordinates
\(X_{e,j}^{(r)}\) with \(e\mid r\).
For such an \(e\), one block has probability generating function
\[
  \alpha_e w^e
  +(1-\alpha_e)
  \left(
    1-\frac1{\delta_r^{(e)}}+
    \frac{w^{eg_r}}{\delta_r^{(e)}}
  \right).
\]
Thus the probability generating function of \(N_r^*\) is
\[
  \Phi_r(w)
  :=
  \prod_{e\mid r}
  \left[
    \alpha_e w^e
    +(1-\alpha_e)
    \left(
      1-\frac1{\delta_r^{(e)}}+
      \frac{w^{eg_r}}{\delta_r^{(e)}}
    \right)
  \right]^{a_e}.
\]
Since \(\Phi_r(w)=\E[w^{N_r^*}]\), for every \(K\ge0\),
\[
  M_{r,K}(q)
  :=
  \lim_{D\to\infty}\E_{F\in U_D}[N_r(F)^K]
  =
  K![t^K]\Phi_r(\exp(t)).
\]
Here moment convergence follows from the common bound
\(0\le N_r(F)\le g_rq^r\).

For \(e\mid r\), put
\[
  \phi_{e,r}(t)
  :=
  \alpha_e\exp(et)
  +(1-\alpha_e)
  \left(
    1-\frac1{\delta_r^{(e)}}+
    \frac{\exp(eg_rt)}{\delta_r^{(e)}}
  \right).
\]
Then \(\Phi_r(\exp(t))=\prod_{e\mid r}\phi_{e,r}(t)^{a_e}\).  Separating the
main term, define
\[
  \mathcal L_r(t)
  :=
  \log\Phi_r(\exp(t))-q^rt
  =
  \sum_{e\mid r}a_e\bigl(\log\phi_{e,r}(t)-et\bigr),
\]
where \(\sum_{e\mid r}ea_e=q^r\).  With \(z=q^rt\),
\[
  M_{r,K}(q)
  =
  q^{rK}K![z^K]\exp(z)
  \exp\left(\mathcal L_r\left(\frac z{q^r}\right)\right).
\]

Until this point \(q\) has been fixed and the limit \(D\to\infty\) has
already been taken.  We now study the resulting moment as \(q\to\infty\).
Fix
\(
  \rho\in(\Z/m\Z)^\times,
\)
and let \(q\to\infty\) through prime powers satisfying
\(q\equiv\rho\pmod m\).  Along such a sequence, the quantities
\(
  g_r,\: \gamma_r^{(e)}\text{and }\delta_r^{(e)}
  \text{ with }e\mid r
\)
are constant.  Indeed, \(g_r\) depends only on \(q^r\bmod m\).  Moreover,
\[
  c_r^{(e)}=1+q^e+\cdots+q^{r-e}
\]
has a residue class modulo \(m\) determined by \(q\bmod m\).  Since
\(g_r\mid m\), this determines \(c_r^{(e)}\bmod g_r\), and hence
\(\gamma_r^{(e)}\) and \(\delta_r^{(e)}\).  We write \([q^{-j}]_\rho\) for coefficient
extraction in the resulting expansion at \(q=\infty\) along this residue
class. The resulting expansion has a remainder of arbitrary fixed order.

\Needspace{6\baselineskip}
\begin{proposition}\label{prop:q-order-expansion}
Fix \(r,K\ge1\), \(J\ge0\), and
\(\rho\in(\Z/m\Z)^\times\).  As \(q\to\infty\) through prime powers with
\(q\equiv\rho\pmod m\),
\[
  \begin{aligned}
    M_{r,K}(q)
    ={}&
    q^{rK}\sum_{j=0}^Jq^{-j}K![z^K]
    \left(
      \exp(z)\,[q^{-j}]_\rho
      \exp\left(\mathcal L_r\left(\frac z{q^r}\right)\right)
    \right) \\
    &+O_{m,\kappa,r,K,J}(q^{rK-J-1}).
  \end{aligned}
\]
The implied constant may be chosen uniformly over the finitely many reduced
residue classes modulo \(m\).
\end{proposition}

\begin{proof}
Fix \(r,K\), and the residue class \(\rho\). Each coefficient through
\(t\)-degree \(K\) in \(\phi_{e,r}(t)\) is rational in \(q\), as is
\(a_e\). Truncating the logarithm and exponential at degree \(K\) shows
that \(M_{r,K}(q)\) is rational in \(q\) on this class.
The bound \(0\le N_r^*/q^r\le m\) implies that
\(q^{-rK}M_{r,K}(q)\) is bounded as \(q\to\infty\). Hence this rational
function has no pole at \(q^{-1}=0\), and its Taylor expansion there has a
remainder of order \(q^{-J-1}\).

For completeness, the mean formula in Theorem~\ref{thm:mean-covariance}
gives \(\E N_r^*=q^r+O_{m,r}(q^{r/2})\), with equality
\(\E N_1^*=q\) when \(r=1\). Independence and
\(0\le X_e^{(r)}\le em\) give
\(\Var(N_r^*)=O_{m,r}(q^r)\). Thus \(N_r^*/q^r\to1\) in \(L^2\),
and boundedness gives convergence of every fixed moment. The constant
term of the normalized moment expansion is therefore \(1\).
The coefficient identity preceding the proposition identifies its Taylor
coefficients with the displayed expression. Multiplying the remainder by
\(q^{rK}\) proves the claim. There are only finitely many reduced residue
classes, so their remainder bounds may be chosen uniformly.
\end{proof}

Equivalently, along the chosen residue class write
\[
  \mathcal L_r\left(\frac z{q^r}\right)
  =
  q^{-1}A_{r,1}^{(\rho)}(z)
  +q^{-2}A_{r,2}^{(\rho)}(z)
  +q^{-3}A_{r,3}^{(\rho)}(z)
  +O(q^{-4})
\]
coefficientwise through the required \(z\)-degree.  Then
\begin{align*}
    M_{r,K}(q)
    ={}&q^{rK}\\
    &+q^{rK-1}K![z^K]\left(\exp(z) A_{r,1}^{(\rho)}(z)\right)\\
    &+q^{rK-2}K![z^K]\left(
      \exp(z)\left(A_{r,2}^{(\rho)}(z)
      +\frac12A_{r,1}^{(\rho)}(z)^2\right)
    \right)\\
    &+q^{rK-3}K![z^K]\left(
      \exp(z)\left(
        A_{r,3}^{(\rho)}(z)
        +A_{r,1}^{(\rho)}(z)A_{r,2}^{(\rho)}(z)
        +\frac16A_{r,1}^{(\rho)}(z)^3
      \right)
    \right)\\
    &+O_{m,\kappa,r,K}(q^{rK-4}).
\end{align*}

\begin{remark}
Combining the finitely many residue classes gives an equivalent expansion over
all admissible prime powers \(q\), with bounded coefficients that are periodic
functions of \(q\bmod m\).  The classwise formulation above is used so that
\([q^{-j}]_\rho\) is ordinary coefficient extraction.
\end{remark}

\subsection{Examples}

\begin{example}[The squarefree hyperelliptic specialization]
\label{ex:squarefree-hyperelliptic}
Assume \(\kappa=m=2\), so \(q\) is odd.  Then
\[
  \alpha_e=\frac1{q^e+1},
  \qquad
  (1-\alpha_e)=\frac{q^e}{q^e+1},
  \qquad
  g_r=2.
\]
Moreover, for \(e\mid r\),
\[
  \gamma_r^{(e)}
  =
  \begin{cases}
    1,& r/e\text{ is odd},\\
    2,& r/e\text{ is even},
  \end{cases}
  \qquad
  \delta_r^{(e)}
  =
  \begin{cases}
    2,& r/e\text{ is odd},\\
    1,& r/e\text{ is even}.
  \end{cases}
\]
Hence Theorem~\ref{thm:mean-covariance} gives
\[
  \lim_{D\to\infty}
  \E_{F\in U_D}[N_r(F)]
  =
  \begin{dcases}
    q^r, & r\text{ odd},\\[4pt]
    q^r+\sum_{e\mid r/2}a_e e\frac{q^e}{q^e+1},
    & r\text{ even}.
  \end{dcases}
\]
Thus even extension degrees have a positive excess term.

The local covariance term simplifies to
\[
  \sigma_{r,s}^{(e)}
  =
  \begin{dcases}
    e^2\frac{q^e}{q^e+1},
    & r/e\text{ and }s/e\text{ are both odd},\\[4pt]
    0,
    & r/e\text{ and }s/e\text{ have opposite parity},\\[4pt]
    e^2\frac{q^e}{(q^e+1)^2},
    & r/e\text{ and }s/e\text{ are both even}.
  \end{dcases}
\]

In particular, if \(r\) and \(s\) are coprime and odd, then only the
degree-\(1\) block contributes, and
\[
  \lim_{D\to\infty}
  \Cov_{F\in U_D}(N_r(F),N_s(F))
  =
  a_1\frac{q}{q+1}
  =
  \frac{q^2}{q+1}.
\]
\end{example}

\begin{remark}\label{rem:even-excess}
In Example~\ref{ex:squarefree-hyperelliptic}, if \(r\) is even and fixed, the
dominant term in the excess comes from \(e=r/2\). Since
\(ea_e=q^e+O_r(q^{e/2})\), one obtains
\[
  \lim_{D\to\infty}
  \E_{F\in U_D}[N_r(F)]
  =q^r+q^{r/2}+O_r(q^{r/4})
  \qquad(r\text{ even}).
\]
This is the large-\(q\) expansion of the affine average obtained after
\(D\to\infty\). The leading excess term \(q^{r/2}\) parallels the phenomenon
studied by Brock and Granville~\cite{BrockGranville}.
Katz~\cite{KatzQuadraticExcess} explains the quadratic excess via the
Frobenius--Schur indicator and investigates analogous excess in extensions
of higher degree.
\end{remark}

\begin{example}[Moments of \(N_1\) in the squarefree case]
Assume \(\kappa=2\).  We use the notation \(M_{r,K}(q)\) and
\(\mathcal L_r\) introduced above.  The
expansions are first obtained on each reduced residue class modulo \(m\).  The
formulas below are uniform over those finitely many classes, so we retain the
\(q\)-dependent notation \(g_1,g_2,\gamma_2^{(1)}\).

First take \(r=1\).
In this case,
\[
  \mathcal L_1\!\left(\frac zq\right)
  =
  q\log\left(
    \frac{q(g_1-1)+g_1\exp(z/q)+q\exp(g_1z/q)}
         {g_1(q+1)}
  \right)-z.
\]
Expanding coefficientwise in \(q^{-1}\), we get
\[
  \begin{aligned}
    \mathcal L_1\!\left(\frac zq\right)
    &={}\frac{g_1-1}{2q}z^2\\[4pt]
    &\quad+\frac1{q^2}
    \left(
      -\frac{g_1-1}{2}z^2
      +\frac{(g_1-1)(g_1-2)}{6}z^3
    \right)
    +O_{m,K}(q^{-3}).
  \end{aligned}
\]
Therefore
\begin{align*}
    &\exp\!\left(\mathcal L_1\!\left(\frac zq\right)\right)\\*
    &\qquad={}
    1+\frac{g_1-1}{2q}z^2\\
    &\qquad\quad+
    \frac1{q^2}
    \left(
      -\frac{g_1-1}{2}z^2
      +
      \frac{(g_1-1)(g_1-2)}{6}z^3
      +
      \frac{(g_1-1)^2}{8}z^4
    \right)\\
    &\qquad\quad+O_{m,K}(q^{-3}).
\end{align*}
Using \(K![z^K]\exp(z)z^j/j!=\binom Kj\), Proposition~\ref{prop:q-order-expansion}
gives, for fixed \(m\) and \(K\ge2\),
\begin{align*}
    M_{1,K}(q)
    ={}&q^K+\binom K2(g_1-1)q^{K-1}\\
    &-(g_1-1)\binom K2q^{K-2}\\
    &+(g_1-1)(g_1-2)\binom K3q^{K-2}\\
    &+3(g_1-1)^2\binom K4q^{K-2}
    +O_{m,K}(q^{K-3}).
\end{align*}

\end{example}

\begin{example}[The cube-free cubic case]
Assume \(m=\kappa=3\).  Since \(\gcd(3,q)=1\), we have
\(q\equiv1\) or \(q\equiv2\pmod 3\).  In this case
$\alpha_e=\frac{q^e+1}{q^{2e}+q^e+1}$.
For \(N_1\), one always has \(\gamma_1^{(1)}=1\).  Hence
\[
  \lim_{D\to\infty}\E_{F\in U_D}[N_1(F)]=q.
\]
Moreover,
\[
  \lim_{D\to\infty}\Var_{F\in U_D}(N_1(F))
  =
  \begin{dcases}
    \frac{2q^3}{q^2+q+1},
    & q\equiv1\pmod3,\\[4pt]
    0,
    & q\equiv2\pmod3.
  \end{dcases}
\]
Thus, when \(q\equiv2\pmod3\), the cube map on \(\F_q^\times\) is bijective and
\(N_1^*=q\) almost surely.

For \(N_2\), we have \(g_2=\gcd(3,q^2-1)=3\), and
\[
  \gamma_2^{(1)}
  =
  \gcd(3,q+1)
  =
  \begin{cases}
    1,& q\equiv1\pmod3,\\
    3,& q\equiv2\pmod3.
  \end{cases}
\]
Therefore
\[
  \lim_{D\to\infty}\E_{F\in U_D}[N_2(F)]
  =
  \begin{dcases}
    q^2,
    & q\equiv1\pmod3,\\[4pt]
    q^2+\frac{2q^3}{q^2+q+1},
    & q\equiv2\pmod3.
  \end{dcases}
\]

The covariance between \(N_1\) and \(N_2\) is
\[
  \lim_{D\to\infty}
  \Cov_{F\in U_D}(N_1(F),N_2(F))
  =
  \begin{dcases}
    \frac{2q^3}{q^2+q+1},
    & q\equiv1\pmod3,\\[4pt]
    0,
    & q\equiv2\pmod3.
  \end{dcases}
\]
Finally,
\[
  \lim_{D\to\infty}\Var_{F\in U_D}(N_2(F))
  =
  \begin{dcases}
    \frac{2q^3}{q^2+q+1}
    +
    \frac{4(q^2-q)q^4}{q^4+q^2+1},
    & q\equiv1\pmod3,\\[4pt]
    \frac{4q^3(q+1)}{(q^2+q+1)^2}
    +
    \frac{4(q^2-q)q^4}{q^4+q^2+1},
    & q\equiv2\pmod3.
  \end{dcases}
\]

\end{example}

\bibliographystyle{amsplain}
\bibliography{references}

\end{document}